# Infinitely Split Nash Equilibrium Problems in Repeated Games

**Jinlu Li**
Department of Mathematics
Shawnee State University
Portsmouth, Ohio 45662
USA

**Abstract**

In this paper, we introduce the concept of infinitely split Nash equilibrium in repeated games in which the profile sets are chain-complete posets. Then by using a fixed point theorem on posets in [8], we prove an existence theorem. As an application, we study the repeated extended Bertrant duopoly model of price competition.

## 1. Introduction and Preliminaries

1.1. Definitions and notations in game theory

In begin of this section, we review some concepts and notations in game theory that are used by many authors. The readers are referred to [1], [6−7], [10−14] for more details. Let $n$ be a positive integer greater than 1. An $n$-person noncooperative strategic game, simply called an $n$-person game, consists of the following elements:

1. the set of $n$ players denoted by $N$ with $|N| = n$;
2. the set of profiles $S_N = \Pi_{i \in N} S_i$, where $S_i$ is the pure strategy set for player $i \in N$;
3. the utility vector mapping $f = \Pi_{i \in N} f_i: S_N \to R^n$, where $f_i$ is the utility (payoff) function of player $i$, for $i \in N$.

This game is denoted by $G(N, S_N, f)$. Throughout this paper, we always assume that, in the products $\Pi_{i \in N} S_i$, $\Pi_{i \in N} f_i$ and $\Pi_{k \in N \setminus \{i\}} S_k$, the players appear in the same sequential orders. As usual, for every $i \in N$, we often denote a profile of pure strategies for player $i$'s opponents by
$$x_{-i} = (x_1, x_2, \ldots, x_{i-1}, x_{i+1}, x_n).$$

The set of profiles of pure strategies for player $i$'s opponents is then denoted by
$$S_{-i} = \Pi_{k \in N \setminus \{i\}} S_k.$$

Hence we may write $x \in S_N$ as

$$x = (x_i, x_{-i}) \text{ with } x_{-i} \in S_{-i}, \text{ for } i \in N. \tag{1}$$

Moreover, for every $x_{-i} \in S_{-i}$, we denote

$$f_i(S_i, x_{-i}) = \{f_i(z_i, x_{-i}): z_i \in S_i\}.$$

From $f = \Pi_{i \in N} f_i$ in the game $G(N, S_N, f)$, for any $x \in S$, we have

$$f(x) = \Pi_{i \in N} f_i(x).$$

One of the most important topics in game theory is the study of Nash equilibrium problems. It has been widely studied by many authors and has been extensively applied to economic theory, business and related industries (see [1], [7], [10−14]). We recall the definition of Nash equilibrium in *n*-person noncooperative strategic games below.

Let $G(N, S_N, f)$ be an *n*-person game. A profile of pure strategies $\hat{x} = (\hat{x}_1, \hat{x}_2, \ldots, \hat{x}_n) \in S_N$ is a Nash equilibrium of this game if and only if, it satisfies

$$f_i(z_i, \hat{x}_{-i}) \leq f_i(\hat{x}_i, \hat{x}_{-i}), \text{ for every } i \in N \text{ and for every } z_i \in S_i, \tag{2}$$

It can be rewritten as

$$f_i(z_i, \hat{x}_{-i}) \leq f_i(\hat{x}_i, \hat{x}_{-i}), \text{ for every } i \in N \text{ and for every } z \in S_N.$$

In an *n*-person game $G(N, S_N, f)$, we define a mapping $F: S_N \times S_N \to R^n$ by

$$F(z, x) = \Pi_{i \in N} f_i(z_i, x_{-i}), \text{ for any } x, z \in S_N. \tag{3}$$

$F(z, x)$ is called the utility vector at profile $x \in S_N$ associated to $z \in S_N$. It is clearly to see

$$F(x, x) = f(x), \text{ any } x \in S_N. \tag{4}$$

Let $\geq^n$ be the component-wise partial order on $R^n$ satisfying that, for $x, z \in S_N$,

$$F(z, x) \leq^n F(x, x) = f(x), \text{ if and only if, } f_i(z_i, x_{-i}) \leq f_i(x_i, x_{-i}), \text{ for all } i \in N. \tag{5}$$

From (2−5), the Nash equilibrium can be rewritten by: a profile $\hat{x} \in S_N$ is a Nash equilibrium of $G(N, S_N, f)$ if and only if,

$$F(z, \hat{x}) \leq^n F(\hat{x}, \hat{x}) = f(\hat{x}), \text{ for all } z \in S_N. \tag{6}$$

1.2. *n*-person dual games

An *n*-person game $G = (N, S_N, f)$ is static. Some games in the real world may not be static. That is, it may not be one-shot nature. It is more realistic for this game to be repeatedly played. The dynamic model of game based on an *n*-person game $G = (N, S_N, f)$ is formulized by the process that this static game is repeated infinite periods (times). It is called an *n*-person repeated game, in which there is a discount factor involved for the utilities (see [10]). The dynamic model for *n*-person repeated games will be studied in section 3

In this paper, we first consider a special model: *n-person dual game*. An *n*-person dual game based on an *n*-person game $G(N, S_N, f)$ is modeled as follows: At first, the players play the game as a static *n*-person noncooperative strategic game. After this game is played first time and before this game is played again, every player always considers the reaction of its competitors to its strategy applied in the first time. To seek the optimization of the player's utilities, the players may make arrangements of strategies to use in the second play. Suppose that this performance is represented by a mapping $A$ on $S_N$. Hence, if $x \in S_N$ is the profile used by the players in the first time, then $Ax \in S_N$ will be the profile used by the players in the second time. This *n*-person dual game is denoted by $G(N, S_N, f, A)^2$.

Then we ask: Is there a Nash equilibrium $\hat{x} \in S_N$ of the game $G(N, S_N, f)$ (first play) such that $A\hat{x} \in S_N$ is also a Nash equilibrium of this game in second play with respect to the translated profiles? It raises the so called split Nash equilibrium problems for dual games.

In [2–5] and [9], multitudinous of iterative algorithms have provided for the approximations of split Nash equilibria for two games. In all results about estimating Nash equilibria in the listed papers, there is a common essential prerequisite: The existence of a Nash equilibrium in the considered problem is assumed. It is indubitable that the existence of solutions for split Nash equilibrium problems is always the crux of the matter for solving these problems.

In [6], the present author proved an existence theorem of split Nash equilibrium problems for related games by using the Fan-KKM Theorem. Since the present author has studied the fixed point theory on posets for several years and has made some applications to Nash equilibrium problems, so, in this paper, we will apply some fixed point theorems on posets to study the solvability of split Nash equilibrium problems for dual games. To this end, the profile sets of games must be equipped with partial orders that may be neither linear spaces, nor topological spaces. The positive aspect of this research is that the utility functions in the considered games are unnecessary to be continuous and the mapping $A$ that defines the split Nash equilibrium problems is unnecessary to be linear.

In section 3, we extend the concept of split Nash equilibrium problems for dual games to infinitely split Nash equilibrium problems for repeated games and prove an existence theorem. As applications, in section 4, we study the existence of infinitely split Nash equilibrium and Nash equilibrium for the repeated extended Bertrant duopoly Model of price competition that is a special repeated game.

**2. Split Nash equilibrium problems in dual games**

2.1. Definitions and notations for split Nash equilibrium problems in dual games

Let $G(N, S_N, f)$ be an *n*-person game. Throughout this paper, unless otherwise stated, we assume that, for every player $i \in N$, his strategy set $S_i$ is nonempty and is equipped with a partial order $\succcurlyeq_i$. That is, for every $i \in N$, player *i*'s strategy set is assumed to be a poset $(S_i, \succcurlyeq_i)$. As the product partially ordered set of $(S_i, \succcurlyeq_i)$'s, the profile set is also a poset $(S_N, \succcurlyeq)$ in which the partial order $\succcurlyeq$ is the component-wise partial order of $\succcurlyeq_i$'s. That is, for $x = (x_1, x_2, \ldots, x_n)$ and $y = (y_1, y_2, \ldots, y_n) \in S_N$, we have that

$$y \succcurlyeq x, \text{ if and only if, } y_i \succcurlyeq_i x_i, \text{ for all } i \in N.$$

For every $i \in N$, $(S_{-i}, \succcurlyeq_{-i})$ is similarly defined to be the product poset of $(S_j, \succcurlyeq_j)$'s, $j \neq i$, in which $\succcurlyeq_{-i}$ is the corresponding component-wise partial order of $\succcurlyeq_j$'s, $j \neq i$.

**Definition 1.** *Let $G(N, S_N, f, A)^2$ be an n-person dual game. The split Nash equilibrium problem associated with this dual game, denoted by $\text{NSNE}(G(N, S_N, f, A)^2)$, is formalized as: to find a profile $\hat{x} \in S_N$ satisfying*

$$f_i(z_i, \hat{x}_{-i}) \leq f_i(\hat{x}_i, \hat{x}_{-i}), \text{ for every } i \in N \text{ and } z \in S, \tag{7}$$

*such that the profile $A\hat{x} \in S_N$ solves the following*

$$f_i((Az)_i, (A\hat{x})_{-i}) \leq f_i((A\hat{x})_i, (A\hat{x})_{-i}), \text{ for every } i \in N \text{ and } z \in S. \tag{8}$$

From (6), a profile $\hat{x} \in S_N$ satisfying (7–8) can be rewritten as:

$$F(z, \hat{x}) \leq^n F(\hat{x}, \hat{x}) = f(\hat{x}), \text{ for all } z \in S_N, \tag{9}$$

and

$$F(Az, A\hat{x}) \leq^n F(A\hat{x}, A\hat{x}) = f(\hat{x}), \text{ for all } z \in S_N. \tag{10}$$

Such a profile $\hat{x}$ in $S_N$ is called a split Nash equilibrium of this split Nash equilibrium problem $\text{NSNE}(G(N, S_N, f, A)^2)$. The set of all split Nash equilibriums is denoted by $\mathcal{S}(G(N, S_N, f, A)^2)$.

When looking at the equilibrium problems (5) and (6) separately, the problem (5) is the classical Nash equilibriums problem of strategic games. When, considering a special case, $A = I$ ($A$ is unnecessary to be linear), that is the identity mapping on $S_N$, $\text{NSNE}(G(N, S_N, f, I)^2)$ reduces to the classical Nash equilibrium problem for the game $G(N, S_N, f)$. In this view, split Nash equilibrium problems for dual games can be considered as the natural extensions of the classical Nash equilibrium problems.

A fixed point theorem on posets is proved in [8]. In this theorem, the underlying space is a chain-complete poset and the considered mapping is just required to satisfy order-increasing upward condition without any continuity condition (As a matter of fact, the underlying space is just equipped with a partial order and it may not have any topological structure). The values of the considered mapping are universally inductive that is a relatively broad concept. Some properties and examples of universally inductive posets have been provided in [8]. We recall this theorem below that will be used in the proof of the main theorems in this paper.

**Fixed Point Theorem A** (Theorem 3.2 in [8]). *Let $(P, \succcurlyeq^P)$ be a chain-complete poset and let $\Gamma: P \to 2^P \setminus \{\emptyset\}$ be a set-valued mapping satisfying the following three conditions:*

    A1. *$\Gamma$ is $\succcurlyeq^P$-increasing upward;*
    A2. *$(\Gamma(x), \succcurlyeq^P)$ is universally inductive, for every $x \in P$;*
    A3. *There is an element $y_*$ in $P$ and $v_* \in \Gamma(y_*)$ with $y_* \preccurlyeq^P v_*$.*

*Let $\mathcal{F}(\Gamma)$ denote the set of fixed points of $\Gamma$. Then*

(i)  $(\mathcal{F}(\Gamma), \succcurlyeq^P)$ is a nonempty inductive poset;

(ii) $(\mathcal{F}(\Gamma) \cap [y_*), \succcurlyeq^P)$ is a nonempty inductive poset; and $\Gamma$ has an $\succcurlyeq^P$-maximal fixed point $x^*$ with $x^* \succcurlyeq^P y_*$.

## 2.2 An existence theorem for split Nash equilibrium in dual games

We need the following concept, order-positive, for mappings from posets to posets. It is an important condition for the mapping $A$ for the existence of split Nash equilibrium in split Nash equilibrium problems.

**Definition 2.** *Let $(X, \succcurlyeq^X)$, $(Y, \succcurlyeq^Y)$ and $(U, \succcurlyeq^U)$ be posets. Let $C, D$ be nonempty subsets of $X$ and $Y$, respectively. A mapping $g: X \times Y \to U$ is said to be order-positive from $X \times Y$ to $U$ whenever, for $x, y \in D$, if $x \preccurlyeq^Y y$, then*

$$g(z, x) \preccurlyeq^U g(t, x) \quad \text{implies} \quad g(z, y) \preccurlyeq^U g(t, y), \text{ for any } z, t \in X. \tag{11}$$

*In particular, if $(U, \succcurlyeq^U) = (R^m, \geq^m)$, where $m$ is a natural number, a mapping $g$ is order-positive from $X \times Y$ to $(R^m, \geq^m)$ whenever, for $x, y \in C$, if $x \preccurlyeq^X y$, then*

$$g(z, x) \leq^m g(t, x) \quad \text{implies} \quad g(z, y) \leq^m g(t, y), \text{ for any } z, t \in X. \tag{12}$$

Let $G(N, S_N, f, A)^2$ be an $n$-person dual game. To prove an existence theorem for split Nash equilibrium problem $\text{NSNE}(G(N, S_N, f, A)^2)$, we need to define a mapping $\pi: S_N \to 2^{S_N}$, for $x \in S_N$, by

$$\pi(x) = \{t \in S_N: F(z, x) \leq^n F(t, x) \text{ and } F(Az, Ax) \leq^n F(At, Ax), \text{ for all } z \in S_N\}. \tag{13}$$

$\pi$ can be equivalently written, for $x \in S_N$, as

$$\pi(x) = \{t \in S_N: f_i(z_i, x_{-i}) \leq f_i(t_i, x_{-i}) \text{ and } f_i((Az)_i, (Ax)_{-i}) \leq f_i((At)_i, (Ax)_{-i}),$$
$$\text{for every } i \in N \text{ and for all } z \in S_N\}.$$

**Observation 1**. In Theorem 1 given below, it is assumed that, for every $x \in S_N$, $\pi(x) \neq \emptyset$. It means that, for any given profile $x \in S_N$ and for every player $i \in N$, when player $i$'s opponents take $x_{-i}$ to play, there exists a strategy $t_i \in S_i$ such that player $i$ will optimize his utility at the profile $(t_i, x_{-i})$. Hence, the condition that $\pi(x)$ is nonempty is a reasonable condition and it should not be too strong.

Now we prove one of the main theorems of this paper.

**Theorem 1**. *Let $G(N, S_N, f, A)^2$ be an $n$-person dual game. Suppose that, for every $i \in N$, $(S_i, \succcurlyeq_i)$ is a nonempty chain-complete poset. Let $(S_N, \succcurlyeq)$ be the product poset of $(S_i, \succcurlyeq_i)$'s equipped with the component-wise partial order $\succcurlyeq$. If $f$ and $A$ satisfy the following conditions:*

  a). *For every $i \in N$, $f_i$ is order-positive from $(S_i, \succcurlyeq_i) \times (S_{-i}, \succcurlyeq_{-i})$ to $(R, \geq)$;*
  b). *For every $x \in S_N$, $\pi(x)$ is a universally inductive subset of $S_N$;*
  c). *The operator $A: S_N \to S_N$ is $\succcurlyeq$-increasing;*
  d). *There are elements $x' \in S_N$ and $u' \in \pi(x')$ satisfying $x' \preccurlyeq u'$,*

*then the dual game* $G(N, S_N, f, A)^2$ *has a split Nash equilibrium. Moreover*

(i) $(\mathcal{S}(G(N, S_N, f, A)^2), \succcurlyeq)$ *is a nonempty inductive poset;*

(ii) $(\mathcal{S}(G(N, S_N, f, A)^2) \cap [x'), \succcurlyeq)$ *is a nonempty inductive poset.*

*Proof.* Since, for every $i \in N$, $(S_i, \succcurlyeq_i)$ is a nonempty chain-complete poset, then the profile set, as a product space of chain-complete spaces $(S_i, \succcurlyeq_i)$'s, $(S_N, \succcurlyeq)$ is a nonempty chain-complete poset, where $\succcurlyeq$ is the component-wise partial order of $\succcurlyeq_i$'s. Define a set-valued mapping $\Gamma: S_N \to 2^{S_N}$ by

$$\Gamma(x) = \pi(x)$$
$$= \{t \in S_N: F(z, x) \leq^n F(t, x) \text{ and } F(Az, Ax) \leq^n F(At, Ax), \text{ for all } z \in S_N\}, \text{ for } x \in S_N. \quad (14)$$

From condition b) in this theorem, it implies that, for every $x \in S_N$, $\Gamma(x) \neq \emptyset$. Hence the mapping $\Gamma: S_N \to 2^{S_N} \setminus \{\emptyset\}$ is a well-defined set-valued mapping with universally inductive values in $S_N$.

Next we show that $\Gamma$ is $\succcurlyeq$-increasing upward. Notice that the partial order $\succcurlyeq$ on $S_N$ is the component-wise partial order of $\succcurlyeq_i$'s on $S_i$'s, respectively. It implies that, for any $x, y \in S_N$, $x \preccurlyeq y$ is equivalent to $x_i \preccurlyeq_i y_i$ and $x_{-i} \preccurlyeq_{-i} y_{-i}$, for every $i \in N$. From condition a), for every $i \in N$, $f_i$ is order-positive from $(S_i, \succcurlyeq_i) \times (S_{-i}, \succcurlyeq_{-i})$ to $(R, \geq)$. From (3), it implies that $F$ is order-positive from $(S_N, \succcurlyeq) \times (S_N, \succcurlyeq)$ to $(R^n, \geq^n)$. Then, for arbitrary $x, y \in S_N$ with $x \preccurlyeq y$, it implies $x_{-i} \preccurlyeq_{-i} y_{-i}$, every $i \in N$. From condition a), we then have

$$f_i(z_i, x_{-i}) \leq f_i(t_i, x_{-i}) \implies f_i(z_i, y_{-i}) \leq f_i(t_i, y_{-i}), \text{ for } z_i, t_i \in S_i.$$

It follows that

$$F(z, x) \leq^n F(t, x) \quad \text{implies} \quad F(z, y) \leq^n F(t, y), \text{ for any } z, t \in S_N. \quad (15)$$

From condition c), the operator $A: S_N \to S_N$ is $\succcurlyeq$-increasing. It implies that if $x \preccurlyeq y$, then, $Ax \preccurlyeq Ay$. From condition a) again, similar to (15), we have

$$F(Az, Ax) \leq^n F(At, Ax) \quad \text{implies} \quad F(Az, Ay) \leq^n F(At, Ay), \text{ for any } z, t \in S_N. \quad (16)$$

(15) and (16) together imply that if $x \preccurlyeq y$, then $\Gamma(x) \subseteq \Gamma(y)$. Hence $\Gamma$ is $\succcurlyeq$-increasing upward.

The elements $x' \in S_N$ and $u' \in \pi(x')$ given in condition d) in this theorem satisfy that $u' \in \Gamma(x')$ such that $x' \preccurlyeq u'$. So $\Gamma$ satisfies all conditions in the Fixed Point Theorem A. It follows that $\mathcal{F}(\Gamma) \neq \emptyset$ and it satisfies the properties (i) and (ii) in Theorem A. From (9) and (10), the definition of $\mathcal{S}(G(N, S_N, f, A)^2)$, and (14), the definition of $\Gamma$, we obtain

$$\mathcal{S}(G(N, S_N, f, A)^2) = \mathcal{F}(\Gamma). \quad (17)$$

By Applying Theorem A and (17), the proof of this theorem is completed immediately. □

2.3. Applications to partially ordered Banach spaces.

In this subsection, we consider a special case of *n*-person dual games in which the strategy set for every player is a nonempty and compact subset of a partially ordered Banach space. This case should be very useful in the applications. In [8], it was proved that every partially ordered compact Hausdorff topological space is both chain-complete and universally inductive, as a consequence of Theorem 1, we have

**Corollary 1**. *Let $G(N, S_N, f, A)^2$ be an n-person dual game. Suppose that, for every $i \in N$, $S_i$ is a nonempty compact subset of a partially ordered Banach space $(B_i, \succcurlyeq_i)$. Let $B_N = \Pi_{i \in N} B_i$ equipped with the component-wise partial order $\succcurlyeq$ of $\succcurlyeq_i$'s. If f and A satisfy the following conditions*:

   a). *For every $i \in N$, $f_i$ is order-positive from $(S_i, \succcurlyeq_i) \times (S_{-i}, \succcurlyeq_{-i})$ to $(R, \geq)$;*
   b). *For every $x \in S_N$, $\pi(x)$ is a nonempty closed subset of $S_N$;*
   c). *The operator $A: S_N \to S_N$ is $\succcurlyeq$-increasing;*
   d). *There are elements $x' \in S_N$ and $u' \in \pi(x')$ satisfying $x' \preccurlyeq u'$,*

*then the dual game $G(N, S_N, f, A)^2$ has a split Nash equilibrium. Moreover, $\mathcal{S}(G(N, S_N, f, A)^2)$ has the properties* (i) *and* (ii) *listed in Theorem* 1.

**Remarks 2**. In Corollary 1, even though the profile set in the dual game $G(N, S_N, f, A)^2$ is a subset of a Banach space, the operator $A: S_N \to S_N$ is not required to be linear. It may be a nonlinear operator.

### 3. Infinitely Split Nash equilibrium problems in repeated games

3.1. Definitions and notations of *n*-person repeated games

Let $G(N, S_N, f)$ be an *n*-person game. Recall that, for every $i \in N$, player *i*'s strategy set is assumed to be a poset $(S_i, \succcurlyeq_i)$. $(S_N, \succcurlyeq)$ is the product poset, where $\succcurlyeq$ is the component-wise partial order of $\succcurlyeq_i$'s naturally equipped on $S_N$.

For every natural number *k*, after the players repeated play the game *k* times and, for each time, the game is played as a static *n*-person simultaneous-move game, before they play this static game again, every player always considers the reaction of its competitors to its strategy applied in the previous time. To optimize their utilities, the players may make arrangements of strategies to use in the next play (the $(k+1)^{\text{th}}$ play). Suppose that the profile of the arranged strategies is represented by the value of a mapping $A_k: S_N \to S_N$, for $k = 1, 2, 3, \ldots$ (Since $(S_i, \succcurlyeq_i)$ is just a poset, it may not be equipped with any algebraic structure. So the linearity of $A_k$ is may not be defined). To summarizing this process, if $x \in S_N$ is the profile used by the players in the first time, then $A_1 x \in S_N$ will be the profile used by the players in the second time; $A_2 A_1 x \in S_N$ will be the profile used by the players in the third time. Hence, for $k = 1, 2, 3, \ldots$, $A_k \ldots A_2 A_1 x \in S_N$ will be the profile used by the players in the $(k+1)^{\text{th}}$ play. For simplicity, we write

$$\Pi_k = A_k A_{k-1} \ldots A_1 A_0, \text{ for } k = 0, 1, 2, \ldots. \tag{18}$$

where $A_0 = I$, that is the identity operator on $S_N$. Then $\Pi_k: S_N \to S_N$ is a single-valued mapping. In particular, if $A_k = A_{k-1} = \ldots = A_1 = A$, then we denote

$$\Pi_k = A^k A_0, \text{ for } k = 1, 2, \ldots.$$

Suppose that the utilities of this game are bounded. That is, there is a number $M > 0$ such that

(M) $\qquad |f_i(x)| \leq M$, for every $i \in N$ and for all $x \in S_N$.

There is a discount factor $0 < \rho < 1$. For every $i \in N$, player $i$'s discounted value of utility at a profile $x \in S_N$ is

$$h_i(x) = \sum_{k=0}^{\infty} \rho^k f_i(\Pi_k x), \qquad (19)$$

Player $i$'s discounted value of utility at the profile $x$ associated with a profile $z \in S_N$ is

$$H_i(z, x) = \sum_{k=0}^{\infty} \rho^k f_i((\Pi_k z)_i, (\Pi_k x)_{-i}). \qquad (20)$$

It implies

$$H_i(x, x) = h_i(x), \text{ for every } i \in N \text{ and for all } x \in S_N.$$

The utility vector with discounted values for this repeated game at the profile $x$ associated with a profile $z \in S_N$ is

$$H(z, x) = \Pi_{i \in N} H_i(z, x) = \Pi_{i \in N} \left( \sum_{k=0}^{\infty} \rho^k f_i((\Pi_k z)_i, (\Pi_k x)_{-i}) \right).$$

Under the boundedness condition (M), for every $i \in N$, both of $h_i$ and $H_i$ are well-defined real valued functions on $S_N$ and $S_N \times S_N$, respectively. Then it forms an $n$-person dynamic model based on an $n$-person game. It is called an $n$-person repeated game based on the $n$-person static game $G(N, S_N, f)$ and is denoted by

$$G(N, S_N, f, A_k)_{k=0}^{\infty}.$$

**Definition 3.** *Let $G(N, S_N, f, A_k)_{k=0}^{\infty}$ be an $n$-person repeated game. A profile $\hat{x} \in S_N$ is called a Nash equilibrium of this repeated game, if the following inequalities are satisfied*

$$H_i(z, \hat{x}) \leq H_i(\hat{x}, \hat{x}), \text{ for every } i \in N \text{ and } z \in S_N. \qquad (21)$$

**Definition 4.** *The infinitely split Nash equilibrium problem associated with the repeated game $G(N, S_N, f, A_k)_{k=0}^{\infty}$, denoted by $\text{NSNE}(G(N, S_N, f, A_k)_{k=0}^{\infty})$, is formalized as: to find a profile $\hat{x} \in S_N$ satisfying*

$$f_i((\Pi_k z)_i, (\Pi_k \hat{x})_{-i}) \leq f_i((\Pi_k \hat{x})_i, (\Pi_k \hat{x})_{-i}), \text{ for every } i \in N \text{ and } z \in S_N, \text{ for } k = 0, 1, 2, \ldots. \qquad (22)$$

*Such a profile $\hat{x} \in S_N$ is called an infinitely split Nash equilibrium. The set of all infinitely split Nash equilibriums of $G(N, S_N, f, A_k)_{k=0}^{\infty}$ is denoted by $\mathcal{S}(G(N, S_N, f, A_k)_{k=0}^{\infty})$.*

**Proposition 1.** *Every infinitely split Nash equilibrium of an $n$-person repeated game is a Nash equilibrium of this repeated game.*

*Proof.* Suppose that, for every $i \in N$, for $k = 0, 1, 2, \ldots$, the following inequality holds

$$f_i((\Pi_k z)_i, (\Pi_k \hat{x})_{-i}) \leq f_i((\Pi_k \hat{x})_i, (\Pi_k \hat{x})_{-i}), \text{ for all } z \in S_N.$$

Since $0 < \lambda < 1$, it implies

$$H_i(z, \hat{x}) = \sum_{k=0}^{\infty} \lambda^k f_i((\Pi_k z)_i, (\Pi_k \hat{x})_{-i}) \leq \sum_{k=0}^{\infty} \rho^k f_i((\Pi_k \hat{x})_i, (\Pi_k \hat{x})_{-i}) = H_i(\hat{x}, \hat{x}).$$

It completes the proof of this proposition. □

Similar to (13) for the definition of the mapping $\pi$, regarding to $G(N, S_N, f, A_k)_{k=0}^{\infty}$, we need to define a mapping $\psi: S_N \to 2^{S_N}$, for $x \in S_N$, by

$$\psi(x) = \{t \in S_N: f_i((\Pi_k z)_i, (\Pi_k x)_{-i}) \leq f_i((\Pi_k t)_i, (\Pi_k x)_{-i}), \text{ for every } i \in N, \text{ all } z_i \in S_i, \text{ all } k = 0, 1, \ldots\}$$

$\psi$ can be rewritten, for $x \in S_N$, as

$$\psi(x) = \{t \in S_N: F(\Pi_k z, \Pi_k x) \leq^n F(\Pi_k t, \Pi_k x), \text{ for } k = 0, 1, 2, \ldots \text{ and for all } z \in S_N\}. \quad (23)$$

**Theorem 2**. *Let $G(N, S_N, f, A_k)_{k=0}^{\infty}$ be an n-person repeated game. Suppose that, for every $i \in N$, $(S_i, \succcurlyeq_i)$ is a nonempty chain-complete poset. Let $(S_N, \succcurlyeq)$ be the product poset of $(S_i, \succcurlyeq_i)$'s equipped with the component-wise partial order $\succcurlyeq$. Suppose that the following conditions are satisfied*:

  a). *For every $i \in N$, $f_i$ is order-positive from $(S_i, \succcurlyeq_i) \times (S_{-i}, \succcurlyeq_{-i})$ to $(R, \geq)$;*
  b). *For every $x \in S_N$, $\psi(x)$ is a nonempty universally inductive subset of $S_N$;*
  c). *For every $k = 1, 2, \ldots, A_k: S_N \to S_N$ is an $\succcurlyeq$-increasing operator;*
  d). *There are elements $x' \in S_N$ and $u' \in \psi(x')$ satisfying $x' \preccurlyeq u'$.*

*Then the repeated game $G(N, S_N, f, A_k)_{k=0}^{\infty}$ has an infinitely split Nash equilibrium. Moreover*

  (i) $(\mathcal{S}(G(N, S_N, f, A_k)_{k=0}^{\infty}), \succcurlyeq)$ *is a nonempty inductive poset;*
  (ii) $(\mathcal{S}(G(N, S_N, f, A_k)_{k=0}^{\infty}) \cap [x'), \succcurlyeq)$ *is a nonempty inductive poset.*

*Proof*. The proof of this theorem is similar to the proof of Theorem 1. As a product space of chain-complete posets $(S_i, \succcurlyeq_i)$'s, the profile set $(S_N, \succcurlyeq)$ is also a nonempty chain-complete poset, where $\succcurlyeq$ is the component-wise partial orders $\succcurlyeq_i$'s. By using (23), we define a set-valued mapping $\Gamma: S_N \to 2^{S_N}$, for $x \in S_N$, by

$$\Gamma(x) = \psi(x) = \{t \in S_N: F(\Pi_k z, \Pi_k x) \leq^n F(\Pi_k t, \Pi_k x), \text{ for } k = 0, 1, 2, \ldots \text{ and for all } z \in S_N\}. \quad (24)$$

From (3) and (24), $\Gamma(x)$ can be rewritten as

$\Gamma(x) = \psi(x)$

$= \{t \in S_N: f_i((\Pi_k z)_i, (\Pi_k x)_{-i}) \leq f_i((\Pi_k t)_i, (\Pi_k x)_{-i}), \text{ for every } i \in N, \text{ all } z_i \in S_i, \text{ all } k = 0, 1, \ldots\}$

From condition b) in this theorem, the mapping $\Gamma: S_N \to 2^{S_N} \setminus \{\emptyset\}$ is a well-defined set-valued mapping with universally inductive values in $S_N$. Next we show that $\Gamma$ is $\succcurlyeq$-increasing upward. From condition a), for every $i \in N$, $f_i$ is order-positive from $(S_i, \succcurlyeq_i) \times (S_{-i}, \succcurlyeq_{-i})$ to $(R, \geq)$. From condition c), it implies that, for every $k = 0, 1, 2, \ldots$ $\Pi_k: S_N \to S_N$ is $\succcurlyeq$-increasing. Then, for arbitrary $x, y \in S_N$ with $x \preccurlyeq y$, similarly to (15) and (16), we can show that

$$F(\Pi_k z, \Pi_k x) \leq^n F(\Pi_k t, \Pi_k x) \implies F(\Pi_k z, \Pi_k y) \leq^n F(\Pi_k t, \Pi_k y), \text{ for any } z, t \in S_N, k = 0, 1, \ldots. \quad (25)$$

(25) implies that if $x \preccurlyeq y$, then $\Gamma(x) \subseteq \Gamma(y)$. Hence $\Gamma$ is $\succcurlyeq$-increasing upward. The elements $x' \in S_N$ and $u' \in \psi(x')$ given in condition d) in this theorem implies that $u' \in \Gamma(x')$ with $x' \preccurlyeq u'$. So $\Gamma$ satisfies all conditions of Fixed Point Theorem A. Rest of the proof is the same to the proof of Theorem 1. □

Similar to Corollary 1, as an application of Theorem 2 to partially ordered Banach spaces, we have

**Corollary 2**. *Let $G(N, S_N, f, A_k)_{k=0}^{\infty}$ be an n-person repeated game. Suppose that, for every $i \in N$, $S_i$ is a nonempty compact subset of a partially ordered Banach space $(B_i, \succcurlyeq_i)$. Let $(B_N, \succcurlyeq)$ be the product partially ordered Banach space of $(B_i, \succcurlyeq_i)$'s, where $\succcurlyeq$ is the component-wise partial order of $\succcurlyeq_i$'s. Suppose that the following conditions are satisfied*:

a). *For every $i \in N$, $f_i$ is order-positive from $(S_i, \succcurlyeq_i) \times (S_{-i}, \succcurlyeq_{-i})$ to $(R, \geq)$*;
b). *For every $x \in S_N$, $\psi(x)$ is a nonempty closed subset of $S_N$*;
c). *For every $k = 1, 2, \ldots$, $A_k: S_N \to S_N$ is an $\succcurlyeq$-increasing operator,*
d). *There are elements $x' \in S_N$ and $u' \in \psi(x')$ satisfying $x' \preccurlyeq u'$.*

*Then the repeated game $G(N, S_N, f, A_k)_{k=0}^{\infty}$ has an infinitely split Nash equilibrium. Moreover $\mathcal{S}(G(N, S_N, f, A_k)_{k=0}^{\infty})$ has the properties* (i) *and* (ii) *given in Theorem* 2.

By using Proposition 1, as applications of Theorem 2, or in Corollary 2, we obtain the following existence results about Nash equilibrium of *n*-person repeated games.

**Corollary 3**. *Let $G(N, S_N, f, A_k)_{k=0}^{\infty}$ be an n-person repeated game as given in Theorem* 2 (*or in Corollary* 2). *If conditions* a)-d) *listed in Theorem* 2 (*or in Corollary* 2) *are satisfied, then this repeated game has a Nash equilibrium.*

**Remarks 2**. Theorems 1, 2 and Corollary 1, 2, 3 provide some conditions for the existence of infinitely split Nash equilibrium or Nash equilibrium in repeated games. Notice that these conditions are just necessary conditions and are not sufficient conditions. Hence, if the conditions of these existence results do not hold for some repeated games, there still may exist an infinitely split Nash equilibrium. It only means that it cannot be assured that there is one, if these conditions are not satisfied.

4. **Applications to repeated extended Bertrant duopoly model of price competition**

In [6], the present author generalized the Bertrant duopoly model of price competition with two firms from the same price model (see [10]) to the model with possibly different prices. Then the dual extended Bertrant model is introduced and an existence theorem of split Nash equilibrium for the Markov dual extended Bertrant duopoly model of price competition is proved in [6]. We review this duopoly model below.

The extended Bertrant duopoly model of price competition is a model of oligopolistic competition that deals with two profit-maximizing firms, named by 1 and 2, in a market. In this model, it is assumed that the two firms have constant returns to scale technologies with costs $c_1 > 0$ and $c_2 > 0$, per unit produced, respectively, where the costs $c_1$ and $c_2$ are possibly different. Without loss of the generality, we assume

$$c_1 \leq c_2, \tag{26}$$

The inequality (26) means that the qualities of the products by these two firms may be different. More precisely, the quality of the products in firm 1 may not be as good as the quality of the products in firm 2.

Let $p_j$ be the price of the products by firm $j$, for $j = 1, 2$. Let $\delta(p_1, p_2)$ be the demand function in this duopoly market. Let $\delta_j(p_1, p_2)$ be the sale function for firm $j$, for $j = 1, 2$. $f$ and $\delta_j$ are assumed to be continuous functions of two variables and strictly decreasing with respect to every given variable. Suppose that there are positive numbers $\bar{p}_j$, for $j = 1, 2$, such that, for all $p_k$

$$\delta_j(p_j, p_k) \geq 0, \text{ for all } p_j \in [0, \bar{p}_j) \quad \text{and} \quad \delta_j(p_j, p_k) = 0, \text{ for all } p_j \geq \bar{p}_j. \tag{27}$$

Suppose that the socially optimal (competitive) output level in this market is strictly positive and finite for every firm

$$0 < \delta(c_1, c_2) < \infty.$$

For given prices $p_1, p_2$, set by firms 1 and 2, respectively, the market is assumed to be clear. That is,

$$\delta(p_1, p_2) = \delta_1(p_1, p_2) + \delta_2(p_1, p_2).$$

Let $\lambda = c_1/c_2$, that defines the ratio of the qualities of the products by firm 1 to firm 2. From the assumption (17), we have $\lambda \in (0, 1]$. Considered as a noncooperative strategic game, the competition takes place as follows: The two firms simultaneously name their prices $p_1, p_2$, respectively. The sales $\delta_1(p_1, p_2)$ and $\delta_2(p_1, p_2)$ are then satisfied

$$\frac{\delta_1(p_1, p_2)}{\delta(p_1, p_2)} = \begin{cases} 0, & \text{if } p_1 > \lambda p_2 \\ \dfrac{c_1}{c_1 + c_2}, & \text{if } p_1 = \lambda p_2 \\ 1, & \text{if } p_1 < \lambda p_2, \end{cases} \tag{28}$$

and

$$\frac{\delta_2(p_1, p_2)}{\delta(p_1, p_2)} = \begin{cases} 1, & \text{if } p_1 > \lambda p_2 \\ \dfrac{c_2}{c_1 + c_2}, & \text{if } p_1 = \lambda p_2 \\ 0, & \text{if } p_1 < \lambda p_2. \end{cases} \tag{29}$$

We assume that the firms produce to order and so they incur production costs only for an output level equal to their actual sales. Therefore, for given prices $p_1$, $p_2$, the firm $j$ has profits

$$u_j(p_1, p_2) = (p_j - c_j)\delta_j(p_1, p_2), \text{ for } j = 1, 2. \tag{30}$$

In [6], an existence theorem for Nash equilibrium of the extended Bertrant duopoly model is proved. We recalled it below for easy reference.

**Theorem 6.1 in [6]**. *In the extended Bertrant duopoly Model, there is a unique Nash equilibrium $(\hat{p}_1, \hat{p}_2)$. In this equilibrium, both firms set their prices equal to their costs, respectively*: $\hat{p}_1 = c_1$, $\hat{p}_2 = c_2$.

This extended Bertrant duopoly model of price competition with two firms is a 2-person static game. It is denoted by G($N$, $S_N$, $u$), where $N = \{1, 2\}$, $S_j \in [0, \bar{p}_j]$, $u = (u_1, u_2)$, and $u_j$ is defined by (30), for $j = 1, 2$, respectively.

For every natural number $k$, after the two firms repeated play the game $k$ times and, for each time, the game is played as a static 2-person simultaneous-move game and before they name their prices to play again, every firm always considers the reaction of its competitor to its strategy (price) applied in the previous time. Suppose that when this game is played in the $k^{\text{th}}$ time, the two firms names their prices as $p_j^{k-1}$, for $j = 1, 2$, respectively. To optimize their utilities, for example, firm 1 could try to increase its profits by increasing the price from $p_1^{k-1}$ to $p_1^k$ (never excess $p_2^{k-1}$), even though decreasing its sales. Meanwhile, firm 2 could try to increase its profits by decreasing the price from $p_2^{k-1}$ to $p_2^k$ (never lower than $p_1^{k-1}$) for increasing its sales. Suppose that such performance is defined by a linear transformation $A_k$ from $\left(p_1^{k-1}, p_2^{k-1}\right)$ to $\left(p_1^k, p_2^k\right)$. Here we assume $p_1^0 = p_1$, and $p_2^0 = p_2$, that are the prices set by the two firms in the very first time. So, for $k = 1, 2, \ldots$, there is a 2×2 matrix $M_k$:

$$M_k = \begin{pmatrix} \alpha_k & 1 - \beta_k \\ 1 - \alpha_k & \beta_k \end{pmatrix}, \tag{31}$$

where $0 \leq \alpha_k$, $\beta_k \leq 1$, such that

$$(p_1^k, p_2^k) = A_k((p_1^{k-1}, p_2^{k-1})) = \begin{pmatrix} \alpha_k & 1 - \beta_k \\ 1 - \alpha_k & \beta_k \end{pmatrix}(p_1^{k-1}, p_2^{k-1}). \tag{32}$$

It implies

$$0 \leq p_1^{k-1} \leq p_1^k \text{ and } 0 \leq p_2^k \leq p_2^{k-1}, \text{ for } k = 1, 2, \ldots. \tag{33}$$

Hence the process of repeatedly playing the static game $G(N, S_N, u)$ with the sequence of linear transformations $\{A_k\}_{k=0}^{\infty}$ defined by (32) is a dynamic game, that is the repeated extended Bertrant duopoly model of price competition. It is a special repeated game denoted by $G(N, S_N, u, A_k)_{k=0}^{\infty}$.

From Definition 4, $\hat{p} = (\hat{p}_1, \hat{p}_1) \in S_N$ is an infinitely split Nash equilibriums of the repeated game $G(N, S_N, u, A_k)_{k=0}^{\infty}$, if it satisfies

$$u_i((\Pi_k p)_i, (\Pi_k \hat{p})_{-i}) \le u_i((\Pi_k \hat{p})_i, (\Pi_k \hat{p})_{-i}), \text{ for every } i = 1, 2, \text{ all } p \in S_N, \text{ every } k = 0, 1, 2, \ldots. \quad (34)$$

**Theorem 3**. *For the infinitely split Nash equilibrium problem of the repeated extended Bertrant duopoly model $G(N, S_N, u, A_k)_{k=0}^{\infty}$, we have*

(i) *If $c_1 = c_2 = c$, then, for any sequence of linear transformations $\{A_k\}_{k=0}^{\infty}$ defined in (32), $\hat{p} = (\hat{p}_1, \hat{p}_1) = (c, c)$ is the unique infinitely split Nash equilibrium;*

(ii) *If $c_1 < c_2$, then there exists a unique infinitely split Nash equilibrium $\hat{p} = (\hat{p}_1, \hat{p}_1) = (c_1, c_2)$, only if all linear transformations $A_k$'s equal to the identity, that is*

$$A_k = \begin{pmatrix} 1 & 0 \\ 0 & 1 \end{pmatrix}, \text{ for } k = 0, 1, 2, \ldots. \quad (35)$$

*Proof*. Part (i) is an immediate consequence of Theorem 6.1 in [6]. To prove part (ii), notice that every infinitely split Nash equilibrium of the repeated game $G(N, S_N, u, A_k)_{k=0}^{\infty}$ is a split Nash equilibrium of the dual game $G(N, S_N, u, A_1)^2$ studied in [6]. From Theorem 6.2 in [6], $(\hat{p}_1, \hat{p}_1) = (c_1, c_2)$ is the unique split Nash equilibrium of the dual game $G(N, S_N, u, A_1)^2$, only if

$$A_1 = \begin{pmatrix} 1 & 0 \\ 0 & 1 \end{pmatrix}.$$

It implies $\Pi_2 = A_2$. Applying Theorem 6.1 in [6] again, it follows that $(\hat{p}_1, \hat{p}_1) = (c_1, c_2)$ is the unique split Nash equilibrium of the triple game $G(N, S_N, u, I, A_2)^3$, only if

$$A_2 = \begin{pmatrix} 1 & 0 \\ 0 & 1 \end{pmatrix}.$$

Then (35) is proved by induction. □

Since $\delta_1$ and $\delta_2$ are continuous functions, from the condition (27) and the definition (30) of the utility functions, it implies that there exists $M > 0$, such that

(M) $\quad\quad\quad |u_i(p)| \le M$, for $i = 1, 2$ and for all $p \in S_N$.

Let $\rho$ be the discount factor of this dynamic game. By Proposition 1, as a consequence of Theorem 3, we have

**Corollary 4**. *In the repeated extended Bertrant duopoly model* $G(N, S_N, u, A_k)_{k=0}^{\infty}$, *there is a Nash equilibrium* $(\hat{p}_1, \hat{p}_1) = (c_1, c_2)$ *at which, every firm has zero discounted value of utility, i.e.*

$$h_i(\hat{p}_1, \hat{p}_1) = \sum_{k=0}^{\infty} \rho^k u_i\left(\Pi_k\left(\hat{p}_1, \hat{p}_1\right)\right) = \sum_{k=0}^{\infty} \rho^k u_i(c_1, c_2) = 0, \text{ for } i = 1, 2.$$

**References**

[1] Bade, S., Nash equilibrium in games with incomplete preferences, *Rationality and Equilibrium Studies in Economic Theory*, Volume 26, (2006) 67–90.

[2] Bnouhachem, A., Strong convergence algorithms for split equilibrium problems and Hierarchical fixed point theorems, *The scientific world journal*, Vol. 2014, Article ID 390956.

[3] Bui, D,, Son, D. X., Jiao, L., and Kim, D. S., Line search algorithms for split equilibrium problems and nonexpansive mappings, *Fixed Point Theory and Applications*, 2016:27 (2016).

[4] Censor, Y., Gibali, A. and Reich, S., Algorithms for split variational inequality problem, *Numerical Algorithms*, **59** (2012) 301–323.

[5] Chang, S. S., Wang, L., Kun, Y., and Wang, G., Moudafis open question and simultaneous iterative algorithm for general split equality variational inclusion problems and general split equality optimization problems, *Fixed Point Theory and Applications*, (2014) 301–215.

[6] Li, J. L., Split Equilibrium Problems for Related Games and Applications to Economic Theory, submitted.

[7] Li, J. L., Several Extensions of the Abian-Brown Fixed Point Theorem and Their Applications to Extended and Generalized Nash Equilibria on Chain-Complete Posets, *Journal of Mathematical Analysis and Applications*, **409**, (2014) 1084–1092.

[8] Li, J. L., Inductive properties of fixed point sets of mappings on posets and on partially ordered topological spaces, *Fixed Point Theory and Applications*, (2015) 2015:211 DOI 10.1186/s13663-015-0461-8.

[9] Ma, Z., Wang, L., Chang, S. S., and Duan, W., Convergence theorems for split equality mixed equilibrium problems with applications, *Fixed Point Theory and Applications*, 2015:31 (2015).

[10] Mas-Colell, A., and Whinston, M. D., *Microeconomics Theory*, Oxford University Press, (1995).

[11] Osborne, M. J., *An introduction to game theory*, Oxford University Press, New York, Oxford, 2004.

[12] Xie, L. S., Li, J. L., and Yang, W. S., Order-clustered fixed point theorems on chain-complete preordered sets and their applications to extended and generalized Nash equilibria, *Fixed Point Theory and Applications*, 2013, **2013**/1/192, 13, 1687–1812.

[13] Xu, H. K., A variable Krasnosel'skii-Mann algorithm and the multiple-sets split feasibility problem, *Inverse Problem*, **22** (2006) 2021–2034.

[14] Zhang, C. J., *Set-Valued Analysis with Applications in Economics*, Sciences Press, Beijing (2004) (in Chinese).